\documentclass[a4paper,12pt]{article}
\usepackage{mathrsfs}
\usepackage{amsmath,amssymb,amsthm,amsfonts}
\begin{document}

\newtheorem{theorem}{Theorem}[section]
 \newtheorem{corollary}[theorem]{Corollary}
 \newtheorem{lemma}[theorem]{Lemma}

\title{\Large \bf On the existence of a rainbow 1-factor in proper coloring of $K_{rn}^{(r)}$
\footnote {Supported by NSFC, PCSIRT and the ``973" program. Email:
 lxl@nankai.edu.cn (Xueliang Li);  irisxuzx@gmail.com (Zhixia Xu). $\dag $ corresponding author} }

\author{\small Xueliang Li$^1$ and Zhixia Xu$^{1,2,\dag }$\\
[2mm]
\small $^1$Center for Combinatorics and LPMC-TJKLC,\\
\small Nankai University, Tianjin 300071, P.R. China\\
\small $^2$College of Mathematics and System Sciences,\\
\small Xinjiang University, Urumuqi, 830046, P.R. China}
\date{}
\maketitle
\begin{abstract}
El-Zanati et al proved that for any 1-factorization $\mathcal{F}$ of
the complete uniform hypergraph $\mathcal {G}=K_{rn}^{(r)}$ with
$r\geq 2$ and $n\geq 3$, there is a rainbow 1-factor. We generalize
their result and show that in any proper coloring of the complete
uniform hypergraph $\mathcal {G}=K_{rn}^{(r)}$ with $r\geq 2$ and
$n\geq 3$, there is a rainbow
1-factor.\\[2mm]
\end{abstract}

\vspace{10pt}

\section{\large Introduction}

A hypergraph $\mathcal {G}=(\mathcal {V},\mathcal {E})$ consists of
a finite set $\mathcal {V}$ of vertices and a set $\mathcal {E}$ of
subsets of $\mathcal {V}$ called edges. An edge subset
$\mathcal{E}'$ of disjoint edges of $\mathcal {E}$ is called {\it
independent}. A proper coloring of $\mathcal {E}$ is a partition of
$\mathcal {E}$ into independent sets with each partition set is
given a color, say $1,2,\cdots,l$. For a given coloring of $\mathcal
{G}$, a subhypergraph $\mathcal{G}'$ is called {\it rainbow} if each
edge of $\mathcal{G}'$ has distinct color.  A 1-factor of a
hypergraph $(\mathcal {V},\mathcal {E})$ is an independent edge set
which partition $\mathcal {V}$. A 1-factorization of $(\mathcal
{V},\mathcal {E})$ is a partition of $\mathcal {E}$ into 1-factors.
For positive integers $r\geq 2$ and $n$, the complete r-uniform
hypergraph on $n$ vertices is the hypergraph $K_n^{(r)}$, with a
vertex set $\mathcal{V}$ of order $n$ and an edge set $\mathcal{E}$
consisting of all $r$-subsets of $\mathcal{V}$. Note that
$K_n^{(2)}$ is $K_n$, the simple complete graph of order $n$. In
order for $K_n^{(r)}$ to contain a 1-factor, it is clearly necessary
that $r$ divides $n$. In 1973, Baranyai \cite{Zs} showed that
$K_{nr}^{(r)}$ has a 1-factorization. Given a 1-factorization
$\mathcal {F}$ of $K_{nr}^{(r)}$, $n\geq 3$, Woolbright in 1978
showed that there exists a 1-factor in $K_{nr}^{(r)}$ whose edges
belong to at least $n-1$ different 1-factors of $\mathcal {F}$. In
1998, Woolbright and Fu \cite{D.E.Woolbright} proved that for any
1-factorization of $K_{2n}$ there is a rainbow 1-factor.

In \cite{S.I.El-Zanati}, El-Zanati et al proved that for any
1-factorization $\mathcal{F}$ of the complete uniform hypergraph
$\mathcal {G}=K_{rn}^{(r)}$ with $r\geq 2$ and $n\geq 3$, there is a
rainbow 1-factor. It is clear that a 1-factorization is a very
special case of proper colorings. In the present paper, we want to
use a weaker condition, the proper coloring condition, to replace
their stronger condition, the 1-factorization condition, and to
generalize their result as follows: for any proper coloring of the
complete uniform hypergraph $\mathcal {G}=K_{rn}^{(r)}$ with $r\geq
2$ and $n\geq 3$, there is a rainbow 1-factor. To show the result,
we divide the proof into three cases: $r=2$ and $n\geq 3,$ $r>2$ and
$n=3,$ $r>2$ and $n>3$. Notice that the proof of Theorem 3 for the
case $r>2$ in \cite{S.I.El-Zanati} can be used directly to show that
for $r>2, n>3$ there is a rainbow 1-factor in any proper coloring of
the complete uniform hypergraph $\mathcal {G}=K_{rn}^{(r)}$. For the
case $r>2, n=3$, we give a prove in Theorem 2.3. The substantial
part of our proof is to show the result for the case $r=2$ and
$n\geq 3$, which will be given in Theorems 2.1 and 2.2. As a result
we have

\begin{theorem} For any proper coloring of the complete
uniform hypergraph $\mathcal {G}=K_{rn}^{(r)}$ with $r\geq 2$ and
$n\geq 3$, there is a rainbow 1-factor.
\end{theorem}

\section{\large Existence of a rainbow 1-factor in
 proper coloring of $K_{2n}$ and $K_{3r}^{(r)}$}

 Let $G=(V,E)$ be a graph and $C$ be a proper coloring of $G$.
  A {\it rainbow matching} of $G$ is a
matching whose edges have pairwise different colors. For $e\in E$,
let $C(e)$ denote the color of $e$. For $v\in V$, let $C(v)=\{ C(e)|
e\ \text{is incident with}\ v\}$. For any subset $E'$ of $E$, let
$C(E')=\{ C(e)| e\in E' \}$ and $F(E')=C(E)-C(E')$. For any
$V'\subseteq V$, let $G[V']$ denote the subgraph induced by $V'$.

\begin{lemma} For any proper coloring of $K_{2n}$, there is a
rainbow perfect matching when $n=3$ or $n=4$.
\end{lemma}

{\flushleft\it Proof.}\quad For $n=3$, let the vertices of $K_6$ be
$v_1,v_2,\cdots, v_6$ and $C$ be a proper coloring of $K_6$, let
$1,2,\cdots,l$ be the colors used, we will show that there is a
rainbow $3K_2$ in $C$. Suppose $1$ is the color that appears least
times in $C$ on $E(K_6)$. If $1$ appears on three edges, then $C$ is
a 1-factorization of $K_6$ and by the proof in [1], there is a
rainbow 1-factor in $C$. If 1 appears on two edges, say
$C(x_1x_2)=C(x_3x_4)=1$, assume that $C(x_5x_1)=2,
C(x_5x_2)=3,C(x_5x_3)=4,C(x_5x_4)=5,C(x_5x_6)=6$. Since
$\{x_5x_1,x_6x_2,x_3x_4\}$ is independent, to avoid the existence of
rainbow perfect matching, it must be that $C(x_6x_2)=2$. Similarly,
$\{x_5x_2,x_6x_1,x_3x_4\}$ is independent and $C(x_6x_1)=3$;
$\{x_1x_2,x_5x_3,x_6x_4\}$ is independent and $C(x_6x_4)=4$;
$\{x_1x_2,x_5x_4,x_6x_3\}$ is independent and $C(x_6x_3)=5$. Now
both $\{x_5x_2,x_6x_3,x_1x_4\}$ and $\{x_5x_3,x_6x_2,x_1x_4\}$ are
independent, whatever color the edge $x_1x_4$ receives, we will have
a rainbow perfect matching. If 1 appears only once in $C$ and
$c(x_1x_2)=1$, to avoid the existence of rainbow perfect matching,
there is no rainbow $2K_2$ in the subgraph induced by $\{
v_2,v_3,v_4,v_5\}$. The only such coloring is
$c(x_3x_4)=c(x_5x_6)=2$, $c(x_3x_5)=c(x_4x_6)=3$,
$c(x_3x_6)=c(x_5x_4)=4$. Assume $c(x_3x_1)=5$. Since both
$\{x_1x_3,x_5x_6,x_2x_4\}$ and $\{x_1x_3,x_4x_6,x_2x_5\}$ are
independent, we have $C(x_2x_4)=5$ and $C(x_2x_5)=5$, a
contradiction. So in any proper coloring of $K_6$, there is a
rainbow perfect matching.

For $n=4$, let the vertices of $K_8$ be $v_1,v_2,\cdots, v_8$ and
$C$ be a proper coloring of $K_8$, we will show that there is a
rainbow $4K_2$ in $C$. Starting with any triangle, it is possible to
find in $C$ of $G=K_8$ at least two rainbow $K_4$, say
$G[\{v_1,v_2,v_3,v_4\}]$ and $G[\{v_1,v_2,v_3,v_5\}]$ are both
rainbow. If there is at least one rainbow $2K_2$ in
$G[\{v_5,v_6,v_7,v_8\}]$, since $\{ v_1v_2,v_3v_4\}$, $\{
v_1v_3,v_2v_4\}$ and $\{ v_1v_4,v_2v_3\}$ are all independent and
each edge has a distinct color, we can find a rainbow $4K_4$ in $C$.
Similarly, there is no rainbow $2K_2$ in $G[\{v_4,v_6,v_7,v_8\}]$.
But it is impossible that both $G[\{v_5,v_6,v_7,v_8\}]$ and
$G[\{v_4,v_6,v_7,v_8\}]$ have no rainbow $2K_2$, and the proof is
complete.
 \qed

\begin{theorem} For $n\geq 3$, any proper coloring of $K_{2n}$
contains a rainbow perfect matching.
\end{theorem}

{\flushleft\it Proof.}\quad By Lemma 1, we can assume $n\geq 5$. Let
$C$ be any proper coloring of $K_{2n}$ with the colors named
$1,2,\cdots, l$, $l\geq 2n-1$. Let $\mathcal {M}$ be any maximal
rainbow matching with $|\mathcal {M}|=k$. Suppose $k<n$, we will
show that there must be a rainbow matching with $k+1$ edges. Recall
that $C(\mathcal {M})$ denotes the set of colors of $\mathcal {M}$
and $F(\mathcal {M})$ denotes the complementary set of colors. Let
$s,t$ be two unmatched vertices. We may assume that $C(\mathcal
{M})=\{1,2,\cdots, k\}$ and $C(st)=1$. Note that by maximality of
$\mathcal {M}$, any edge incident with $s$ whose color is in
$F(\mathcal {M})$ must be incident with an edge of $\mathcal {M}$.

Let $C(s)=C_1\cup C_2 $ with $C_1 \subseteq C(\mathcal {M})$, $C_2
\subseteq F(\mathcal {M})$, $|C_1|=p\leq k$, $|C_2|=2n-1-p$.
Consider all the $(s,t)$-paths of length three, whose first edge is
colored with a color $\alpha$ in $F(\mathcal {M})$, and the second
edge is in $\mathcal {M}$; we call them the candidate 3-paths
relative to $\mathcal {M}$. We can assume that each of these paths
has its third edge colored with a color either in $C(\mathcal
{M})-\{1\}$, or the color $\alpha$ again, for otherwise we could
augment $\mathcal {M}$ to $k+1$ edges simply by deleting the second
edge of the path from it and adding the first and third edges. There
are $2n-1-p$ of these candidate paths and only $k-1$ colors in
$C(\mathcal {M})-\{1\}$. So it follows that at least $2n-p-k$ of
these paths have the first and third edges colored with the same
color in $F(\mathcal {M})$; we call such paths $\mathcal
{M}$-symmetric $(s,t)$-paths.

Let $C(t)=C_1'\cup C_2' $ with $C_1' \subseteq C(\mathcal {M})$,
$C_2' \subseteq F(\mathcal {M})$, $|C_1'|=q\leq k$, $|C_2'|=2n-1-q$.
Consider the $2n-1-q$ edges incident with $t$ whose colors are in
$F(\mathcal {M})$. Each of these edges must be incident with an edge
of $\mathcal {M}$, by the maximality of $\mathcal {M}$; at most 2 of
them, say the ones colored $k+1$ and $k+2$, are incident with the
edge of $\mathcal {M}$ colored 1. Now let $L=C(t)\backslash
(C(\mathcal {M})\cup \{k+1,k+2\}$ and $|L|=2n-q-3$.

For each color $i\in L$, we define a slight variation of the
$(\mathcal {M},st)$ pair. If the edge of color $i$ incident with
vertex $t$ is $e_t=\{t,z_i\}$ of $\mathcal {M}$, we let the
corresponding matching be $\mathcal {M}_i=(\mathcal {M}-\{e_i\})\cup
\{e_t\}$; now $t_i$ is unmatched (in $\mathcal {M}_i$), and we let
our starting/ending vertex pair be $s,t_i$, respectively. Note that
$F(\mathcal {M}_i)=(F(\mathcal {M})-{i})\cup \{C(e_i)\}$. Also note
that $C(e_i)\neq 1$, because $i$ is neither $k+1$ or $k+2$.

As in the previous discussion, for each such $i$, there are $2n-1-p$
candidate 3-paths relative to $\mathcal {M}_i$ starting at $s$,
ending at $t_i$, whose first edge is colored with a color in
$F(\mathcal {M}_i)$, and whose second edge is in $\mathcal {M}_i$.
Again, we assume that at least $2n-p-k$ of these paths are
symmetric. Thus, listing the symmetric paths for $i\in L$, we get a
total of at least $(2n-q-3)(2n-p-k)$ paths in the list of symmetric
candidate paths. However, because in each of these symmetric paths
either the middle edge is in $\mathcal {M}$, or the path has the
form $sz_itt_i$, and therefore has the same first and third edges as
$sz_it_it$, each of the symmetric candidate paths is uniquely
determined by its first edge. Moreover, the color $\alpha$ of the
starting/ending edge in these symmetric paths cannot be
$c(e_i)=c(\{z_i,t_i\})$ (the only possible such path has vertex
sequence $stz_it_i$, which is not symmetric because $c(e_i)\neq 1$),
so $\alpha$ must be in $C_2$. Therefore each of the possible
starting color can only start one path in the list. It follows that
$2n-p-1\geq (2n-q-3)(2n-p-k)$. Let $x=2n-p$, $y=2n-q$, then $x-1\geq
(y-3)(x-k)$, and $x(y-4)< k(y-3)$. Since $q\leq k <n$,
$y=2n-q=n+n-q\geq n+n-k >n $, $y>4$. Then we have $\frac{x}{k} <
\frac{y-3}{y-4}$ and $x<k$, that is $p+k>2n$, which is contrary to
$p\leq k <n$.

We conclude that there must be a rainbow matching with $k+1$ edges,
and so the result follows. \qed

For the case $n=3$, it is easy to see that in any 1-factorization of
$K_{3r}^{(r)}$ there is a rainbow 1-factor, and the proof was
omitted in \cite{S.I.El-Zanati}. But in a proper coloring, the proof
is not straightforward, and we prefer to give the details in the
following.

\begin{theorem} For $r\geq 2$, any proper coloring of $K_{3r}^{(r)}$
contains a rainbow perfect matching.
\end{theorem}

{\flushleft\it Proof.}\quad Let the vertex set of $K_{3r}^{(r)}$ be
$\mathcal {V}= \{x_1,x_2,\cdots,x_{3r}\}$ and $C$ be a proper
coloring of $K_{3r}^{(r)}$. Take any two independent edges having a
same color 1, say $m_1= \{x_1,x_2,$ $\cdots,x_r\}$ and $m_2=
\{x_{r+1},x_{r+2},\cdots,x_{2r}\}$. Then for any edge $m^1 \subset
\mathcal {V} - m_1$ other than $m_2$ and $m^{1*}= \mathcal
{V}-(m_1\cup m^1)$, $m^1$ and $m^{1*}$ have the same color and there
is no other edge in this color, otherwise $\{m_1,m^1,m^{1*}\}$ is a
rainbow 1-factor. Similarly, for any edge $m^{2}\subset \mathcal {V}
- m_2$ other than $m_1$ and $m^{2*}=\mathcal {V} -(m_2\cup m^2)$,
$m^2$ and $m^{2*}$ have the same color and there is no other edge in
this color. Let $ m^1=\{x_{r+1},x_{r+2},\cdots,x_{2r-1},x_{2r+1}\}$,
$m^2=\{x_1,x_{2r+2},x_{2r+3},\cdots,x_{3r}\}$, then
$\{m^1,m^2,\{x_2,x_3,\cdots,x_r,x_{2r}\} \}$ is a rainbow 1-factor.
\qed

\end{document}